\documentclass[reqno,11pt,a4paper]{amsart}

\usepackage[utf8]{inputenc}
\usepackage[T1]{fontenc}
\usepackage{multirow}
\usepackage{ae}
\usepackage{units}
\usepackage{icomma}
\usepackage{color}
\usepackage{graphicx}
\usepackage{cancel}
\usepackage{fullpage}
\usepackage{upgreek}
\usepackage{algorithm}
\usepackage{algorithmic}
\usepackage{type1cm}
\usepackage{eso-pic}
\usepackage[breaklinks,pdfpagelabels=false]{hyperref}
\usepackage{verbatim}
\usepackage{tikz}
\usetikzlibrary{arrows}

\usepackage{caption}
\newcommand{\bm}[1]{\mbox{\boldmath $ {#1} $}}

\renewcommand{\epsilon}{\varepsilon}

\newtheorem{theorem}{Theorem}[section]
\newtheorem{lemma}[theorem]{Lemma}

\newtheorem{remark}[theorem]{Remark}

\theoremstyle{definition}

\numberwithin{equation}{section}
\numberwithin{theorem}{section}

\makeatletter

\def\@settitle{\begin{center}\baselineskip14\p@\relax\bfseries\Large\@title\end{center}}

\makeatother

\begin{document}

\title[Non-convergence in heterogeneous Hegselmann-Krause dynamics: \\ The bounded influence model] 
      {Non-convergence in heterogeneous Hegselmann-Krause dynamics: \\ The bounded influence model}

      \keywords{Hegselmann-Krause model, heterogeneous influence bounds, periodic orbits}

      \subjclass[2020]{39A60, 93A14, 91D10, 91D30}

      \date{\today}


      \author{Peter Hegarty, Damiano Ognissanti and Edvin Wedin} \address{Department of Mathematical Sciences, Chalmers University of Technology and University of Gothenburg, 41296 Gothenburg, Sweden} 
\email{hegarty@chalmers.se, damiano@chalmers.se, wedin@fripost.org}






\begin{abstract}
  In a recent preprint, we found counterexamples to a number of auxiliary
  conjectures in a 2012 paper of Mirtabatabaei and Bullo, concerning
  convergence in two fundamental models of heterogeneous Hegselmann-Krause
  dynamics, the synchronized bounded confidence (SBC) and synchronized bounded
  influence (SBI) models. Building on our work, Etesami and Touri presented an
  example of an opinion vector which doesn't converge under SBC dynamics,
  resolving the main conjecture from the 2012 paper in that case. Here we
  complete this project by presenting an example of an opinion vector which
  doesn't converge under SBI dynamics. The example was found with the help of AI. 
\end{abstract}


\maketitle

\setcounter{equation}{0}

\setcounter{equation}{0}
\setcounter{section}{0}

$\;$ \\
$\;$ \par

\section{{\bf Introduction}}\label{sect:intro}

In the simplest, scalar version of Hegselmann-Krause opinion dynamics, one has a finite set of agents, indexed by the integers $1,\,2,\,\dots,\,n$ for some $n \in \mathbb{N}$, whose initial opinions are modelled as real numbers $x_i(0)$. Each agent $i$ also comes equipped with a fixed positive real number $r_i$, called its confidence or influence bound, depending on the context.
In the synchronized bounded confidence (SBC) model, agents update their opinions synchronously in discrete time according to the formula
\begin{eqnarray}\label{eq:dyn}
  x_i (t+1) = \frac{1}{|\mathcal{N}_i (t)|} \sum_{j \in \mathcal{N}_i(t)} x_j (t), 
\end{eqnarray}
where the so-called \emph{neighborhood} of agent $i$ at time $t$ is defined as
\begin{eqnarray}\label{eq:nbdSBC}
  \mathcal{N}_{i}(t) =\{j: |x_j (t) - x_i(t)| \leq r_i\}.
\end{eqnarray}
In words, each agent updates its opinion to the average of those in whom it has
confidence at time $t$.

By contrast, in the synchronized bounded influence (SBI) model, the update rule
is as in (\ref{eq:dyn}), but neighborhoods are defined as
\begin{eqnarray}\label{eq:nbdSBI}
  \mathcal{N}_{i}(t) =\{j: |x_j (t) - x_i(t)| \leq r_j\}.
\end{eqnarray}
In this case, each agent updates its opinion to the average of those who
influence it at time $t$.

The original Hegselmann-Krause model \cite{HK} is the homogeneous setting where
all radii $r_i$ are equal, in which the SBC and SBI models coincide. The
heterogeneous SBC model appears to have been first introduced by Lorenz
\cite{L}, and the SBI model by Mirtabatabaei and Bullo \cite{MB}. The latter
paper posed what came to be one of the major open problems in the field of
opinion dynamics. Conjecture 2.1 of \cite{MB} asserts that, for any
initial opinion vector $\bm{x}(0) \in \mathbb{R}^n$, the vectors $\bm{x}(t)$
converge pointwise as
$t \rightarrow \infty$.

The paper \cite{MB} reads as a kind of ``roadmap'' towards attacking this
conjecture, and contains some partial results and a number of auxiliary
conjectures. In a recent paper \cite{HOW}, we presented counterexamples to
some of these auxiliary conjectures and even to one of the theorems in
\cite{MB}. A couple of weeks later, Etesami and Touri \cite{ET} presented a
counterexample to Conjecture 2.1 of \cite{MB}, for the SBC model. Their example
is remarkably simple and elegant. It consists of only $6$ agents, whose initial
opinions are symmetrically placed about a midpoint - the update rule
(\ref{eq:dyn}) guarantees that this symmetry persists over time. They show that
there exists such a configuration $\bm{x}(0)$ for which the opinion vector
evolves periodically with period 10, i.e.: $\bm{x}(t+10) = \bm{x}(t)$ for all
$t \geq 0$. While the precise coordinates of $\bm{x}(0)$ are complicated
rationals (or very large integers after renormalization), any configuration
$\bm{x}'(0)$ that starts off sufficiently close to $\bm{x}(0)$ will behave
similarly, with $||\bm{x}'(t) - \bm{x}(t)||_{\infty} \rightarrow 0$ as
$t \rightarrow \infty$, thus allowing them to write down an explicit
counterexample by hand.  

Their paper doesn't mention the SBI model at all, however, and at the time of
writing we are not aware of any counterexample to the Mirtabatabaei-Bullo
conjecture for SBI having already been publicized. Indeed, our searches for SBI
examples with few agents have turned up empty-handed, which could explain the
absence of immediate further progress.

Here instead we exploit a property of the SBI dynamics which is not shared by
SBC. Two agents $i$ and $j$ are said to \emph{agree} at time $t$ if
$x_i(t) = x_j(t)$. It is easy to see that, in SBI, if two agents agree at time
$t$ then they will continue to do so forever, whereas this need not be the case
in SBC. An opinion \emph{cluster} is a maximal set of agents that agree and the
number of agents in a cluster is called its \emph{weight}. Our idea was to look
for configurations initially consisting of only a few clusters, but with
varying weights. Below, we present an example of a configuration consisting of
$53$ agents, grouped into $5$ clusters, which does not converge under SBI
dynamics.
As in \cite{ET}, there is such a configuration whose evolution is precisely
periodic, with period $55$, but with complicated initial data. So instead we
give
an example with explicit initial data whose behaviour is asymptotically the
same. It is also worth noting (see Remark \ref{rem:perturb} below) that a
small enough perturbation of both the initial opinion vector and the vector of
influence radii does not affect the asymptotic behaviour, so there exist
non-convergent configurations initially consisting of $53$ individual agents,
all with different opinions and different influence radii.
\\
\par The rest of the paper is organized as follows. Section
\ref{sect:example} presents our counterexample to Conjecture 2.1 of \cite{MB}
for the SBI model. Section \ref{sect:future} briefly discusses possible avenues
for future work in this area, in particular for refinement of the work
contained in this and other recent articles (\cite{HOW, ET}).
Finally, Section \ref{sect:AI} discusses the role of AI in the present work.

\vspace{3mm}

\section{{\bf A non-convergent configuration for the SBI dynamics}}\label{sect:example}

Our example consists of $53$ agents, initially grouped into $5$ clusters denoted
$L, \, C, \, V, R_1$ and $R_2$, whose weights and initial opinions are given in
the table below, along with the respective influence radii, which are also
assumed to be the same for all agents in the same cluster. Note that the
numerical values of the opinions and radii are chosen for convenience: one
could renormalize so that initial opinions are in the interval $[0,\,180]$ and
all influence radii are integers. There is also some flexibility in the choice
of $c(0)$ and $v(0)$, as we will show. All approximations in what follows are
given to four decimal places.  

\begin{table}[ht!]
\begin{center}
\begin{tabular}{|c|c|c|c|} \hline
           ${\hbox{Cluster}}$ & ${\hbox{Weight}}$ & ${\hbox{Radius}}$ & ${\hbox{Initial opinion}}$ \\ \hline \hline
$L$ & $1$ & $9$ & $-8.1$ \\ \hline 
$C$ & $32$ & $5.3$ & $c(0) = -\frac{232}{185} \approx -1.2541$ \\ \hline 
$V$ & $4$ & $3$ & $v(0) = \frac{607}{460} \approx 1.3196$ \\ \hline 
$R_1$ & $6$ & $7.2$ & $9.9$ \\ \hline 
$R_2$ & $10$ & $9$ & $9.9$ \\ \hline 
  \end{tabular}
\end{center}
\caption{A configuration consisting of $53$ agents, grouped into $5$ clusters, which does not converge under the SBI dynamics. The stabilizing clusters $L, \, R_1$ and $R_2$ will be shown to never move, and the middle clusters $C$ and $V$ converge to a periodic trajectory with period $55$.}
\label{tab:one}
\end{table}

\newpage

Define the function $f: \mathbb{R} \rightarrow \mathbb{R}^2$ by
\begin{eqnarray}\label{eq:f}
  f(s) = \left( \frac{s-8.1}{37}, \, \frac{s+99}{46} \right)
\end{eqnarray}
and note that $(c(0), \, v(0)) = f(-38.3)$. For any $s$, let $\bm{x}_s \in \mathbb{R}^{53}$ denote the same configuration as in Table \ref{tab:one}, except that $(c(0), \, v(0)) = f(s)$. The rest of this section is devoted to proving the following:   

\begin{theorem}\label{thm:main}
  If $\bm{x}(0) = \bm{x}_s$ for any $s \in I = [-38.5, \, -38.1]$, and $\bm{x}(t) = \bm{x}_s (t)$ evolves according to the SBI dynamics (\ref{eq:dyn}) and (\ref{eq:nbdSBI}), then the following hold:
  \par (a) the clusters $L, \, R_1$ and $R_2$ never move
  \par (b) $(c(t), \, v(t))$ converges to a periodic orbit of period $55$. In particular, the configuration $\bm{x}_{s}(t)$ does not converge as $t \rightarrow \infty$.  
  \end{theorem}

We begin by defining a \emph{box} $B$ by
\begin{eqnarray}\label{eq:box}
  B = \{(c,\,v) \in \mathbb{R}^2 : -1.4 < c < 0.7, \; 1 < v < 4, \; c < v, \; v-c \leq 5.3\}.
\end{eqnarray}
Note that $(c(0), \, v(0)) = f(s) \in B$, for any $s \in I$. We say that agent $i$ \emph{hears} agent $j$ at time $t$ if $|x_i(t) - x_j(t)| \leq r_j$. The first step is

\begin{lemma}\label{lem:boxedin}
  Suppose that at some time $t \geq 0$, none of the clusters $L, \, R_1, \, R_2$ have yet moved or changed weight, and that $(c(t), \, v(t)) \in B$. Then the following all hold at time $t$:
  \par (a) agents in $L$ only hear one another, and so this cluster retains its position and weight even at time $t+1$,
  \par (b) agents in $R_1 \cup R_2$ all hear one another and nobody else, so each of these also retains its position and weight at time $t+1$,
  \par (c) agents in $C$ hear those in $L$, but do not hear those in $R_1 \cup R_2$,
  \par (d) agents in $V$ hear those in both $C$ and $R_2$, but not those in $L$.
\end{lemma}

All of the above are trivial to verify. They imply that, under the assumptions of the lemma, there are only two sources of uncertainty at time $t$:
\par (i) whether or not agents in $V$ hear those in $R_1$. This is the case if and only if $v(t) \geq 2.7$
\par
\begin{figure}[t]
\centering
\includegraphics[width=\linewidth]{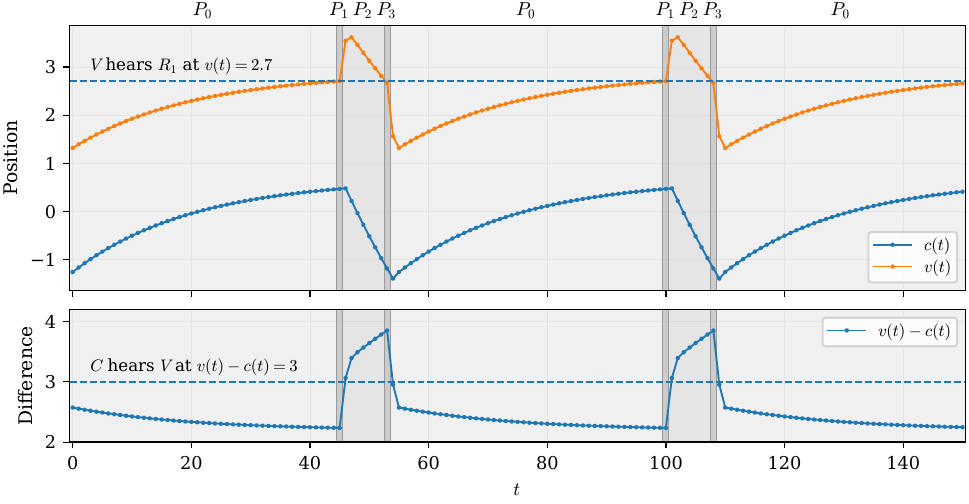}
\caption{The first 150 updates of clusters $C$ and $V$, for the configuration in Table \ref{tab:one}, with the progression through phases $P_0$--$P_3$ highlighted.}
\label{fig:evolution}
\end{figure}

 (ii) whether or not agents in $C$ hear those in $V$. This is the case if and only if $v(t) - c(t) \leq 3$.
\\
\\
Let
\begin{eqnarray}\label{eq:indicators}
  a = \mathbf{1}_{v(t) \geq 2.7}, \;\;\;\; b = \mathbf{1}_{v(t)-c(t) \leq 3}.
\end{eqnarray}
So we identify four possible states at time $t$: the assumptions of Lemma \ref{lem:boxedin} hold and $(a,\,b) \in \{(0,\,1), \, (1,\,1), \, (1,\,0), \, (0,\,0)\}$. These four states are denoted, respectively, as $P_0, \, P_1, \, P_2, \, P_3$, and the progression through them is illustrated in Figure \ref{fig:evolution}. Note that $\bm{x}_s$ is in state $P_0$, for any $s \in I$. Then, 
\begin{eqnarray}\label{eq:boxupdate}
  (c(t+1), \, v(t+1)) = \left\{ \begin{array}{lr} \left( \frac{32c(t) + 4v(t) -8.1}{37}, \, \frac{32c(t) + 4v(t) + 99}{46} \right), & {\hbox{if $\bm{x}(t)$ is in state $P_0$}}, \\ \left( \frac{32c(t) + 4v(t) -8.1}{37}, \, \frac{32c(t) + 4v(t) + 158.4}{52} \right), & {\hbox{if $\bm{x}(t)$ is in state $P_1$}}, \\ \left( \frac{32c(t) -8.1}{33}, \, \frac{32c(t) + 4v(t) + 158.4}{52} \right), & {\hbox{if $\bm{x}(t)$ is in state $P_2$}}, \\ \left( \frac{32c(t) -8.1}{33}, \, \frac{32c(t) + 4v(t) + 99}{46} \right), & {\hbox{if $\bm{x}(t)$ is in state $P_3$}}. \end{array} \right.
\end{eqnarray}

Observe that, if we are in state $P_0$ at time $t$, then 
\begin{eqnarray}\label{eq:onestep}
(c(t+1), \, v(t+1)) = f(z), \;\;\; {\hbox{where $z = 32c(t) + 4v(t)$}}.
\end{eqnarray}
Thus, if $(c(t), \, v(t)) = f(p)$, then $(c(t+1), \, v(t+1)) = f(\lambda p + \mu)$, where $\lambda = \frac{32}{37} + \frac{4}{46} = \frac{810}{851}$
and $\mu = \frac{4 \cdot 99}{46} - \frac{32(8.1)}{37} = \frac{6822}{4255}$. 

The fixed point of the map $p \mapsto \lambda p + \mu$ is $z^{*} = \frac{\mu}{1-\lambda} = \frac{6822}{205} \approx 33.2780$. Starting from 
$p_0 = s$, the $t$:th iterate of this map takes us to 
\begin{eqnarray}\label{eq:precur}
p_t = z^{*} - \lambda^{t} (z^{*} - s).
\end{eqnarray}

Note that the sequence $p_t$ is monotonically increasing whenever $s < z^{*}$. Moreover, one readily verifies that, for all $p \in [-38.5, \, z^{*}]$, 
\begin{eqnarray}\label{eq:recurineq}
-1.4 < \frac{p-8.1}{37} < 0.7, \;\; 1 < \frac{p+99}{46} < 4, \;\; \frac{p+99}{46} - \frac{p-8.1}{37} < 3.
\end{eqnarray}
As we noted above, for any $s \in I$, the initial configuration $\bm{x}(0) = \bm{x}_s$ is in state $P_0$. The inequalities (\ref{eq:recurineq}) imply that the same will be true for $\bm{x}(t)$ as long as $\frac{p_{t} + 99}{46} < 2.7$, that is, $p_{t} < 25.2$, where $p_t$ is given by (\ref{eq:precur}). It is readily checked that, for all $s \in I$, the smallest value of $t$ for which $p_t > 25.2$ is $t=45$. Thus, $\bm{x}(t)$ remains in state $P_0$ for $0 \leq t \leq 44$.
\\
\par At $t=45$ we will have $v(45) \geq 2.7$ and so will enter state $P_1$.
Then, for the next update we will have
\begin{eqnarray*}
(c(46), \, v(46)) = \left(  \frac{p_{46} - 8.1}{37}, \, \frac{p_{46}+158.4}{52} \right).
\end{eqnarray*}
Using (\ref{eq:precur}) one readily checks that
\begin{eqnarray*}
s=p_0 = -38.5 \; \Rightarrow \; c(46) \approx 0.4803, \;\; v(46) \approx 3.5437 \; \Rightarrow v(46)-c(46) > 3, \\
s=p_0 = -38.1 \; \Rightarrow \; c(46) \approx 0.4815, \;\; v(46) \approx 3.5445 \; \Rightarrow v(46)-c(46) > 3.
\end{eqnarray*}
In both cases, we will thus be in state $P_2$ at $t=46$. Since, for any $t$, the configuration $\bm{x}_{s}(t)$ is an affine function of $s$, the same will be true for any $s \in I$. 
\\
\par Moving on, we claim that, for any $s \in I$, we will remain in state $P_2$ for six more time steps, then enter state $P_3$ at $t=53$ before reverting to state $P_0$ at $t=54$. Once again, affineness implies that it suffices to show this when $s=-38.5$ and $s=-38.1$. The next table shows the approximate values of $c(t)$ and $v(t)$ for each $t \in [47, \, 54]$, and both these values of $s$. The calculation can be checked using (\ref{eq:boxupdate}). One merely needs to verify that, in all cases,
\nopagebreak[4]
\begin{itemize}
\item
$(c(t), \, v(t))$ remains inside the box $B$ for all $t$
\item
$v(t) > 2.7$ up to $t=52$, but then $v(53) < 2.7$ and $v(54) < 2.7$
\item
$v(t) - c(t) > 3$ up to $t=53$, but then $v(54) - c(54) < 3$. 
\end{itemize}

\begin{table}[ht!]
  \begin{center}
  \begin{tabular}{|c|c|c|c|c|c|c|c|} \hline
\multirow{2}{*}{Time $t$} & \multicolumn{3}{c|}{$s=-38.5$} & \multicolumn{3}{c|}{$s=-38.1$} & \multirow{2}{*}{State} \\
\cline{2-7}
& $c(t)$ & $v(t)$ & $v(t)-c(t)$ & $c(t)$ & $v(t)$ & $v(t)-c(t)$ & \\
\hline
47 & 0.2203 & 3.6143 & 3.3940 & 0.2214 & 3.6151 & 3.3937 & $P_2$ \\ \hline
48 & -0.0318 & 3.4598 & 3.4916 & -0.0307 & 3.4605 & 3.4912 & $P_2$ \\ \hline
49 & -0.2763 & 3.2927 & 3.5690 & -0.2753 & 3.2934 & 3.5687 & $P_2$ \\ \hline
50 & -0.5134 & 3.1294 & 3.6428 & -0.5124 & 3.1301 & 3.6425 & $P_2$ \\ \hline
51 & -0.7433 & 2.9710 & 3.7142 & -0.7423 & 2.9716 & 3.7139 & $P_2$ \\ \hline
52 & -0.9662 & 2.8173 & 3.7835 & -0.9653 & 2.8179 & 3.7832 & $P_2$ \\ \hline
53 & -1.1824 & 2.6683 & 3.8507 & -1.1815 & 2.6689 & 3.8504 & $P_3$ \\ \hline
54 & -1.3920 & 1.5617 & 2.9537 & -1.3911 & 1.5624 & 2.9535 & $P_0$ \\ \hline
\end{tabular}
\end{center}
\caption{The evolution from $t=47$ to $t=54$ for the mobile clusters $C$ and $V$, in configurations initially satisfying the hypothesis of Theorem \ref{thm:main}.}
\label{tab:two}
\end{table}

Since we are back in state $P_0$ at $t=54$, we can apply (\ref{eq:onestep}) and will have $(c(55), \, v(55)) = f(z)$, where $z=32 \cdot c(54) + 4 \cdot v(54)$. Using the numbers in the last row of Table \ref{tab:two} and affineness, one readily checks that, for any $s \in I$, $z \in (-38.3, \, -38.26) \subset I$. 
\\
\par The proof of Theorem \ref{thm:main} is almost complete. We have shown that, for any $s \in I$, if $\bm{x}(0) = \bm{x}_s$ then $\bm{x}(55) = \bm{x}_{g(s)}$ for some $g(s) \in g(I) \subset [-38.3, \, -38.26] \subset I$. The map $g: I \rightarrow I$ is affine, with slope $\alpha = \frac{g(-38.1)-g(-38.5)}{0.4} \approx 0.0766 < 1$, so $g$ contracts towards a fixed point $s^{*} \approx -38.2805$. If $\bm{x}(0) = \bm{x}_{s^{*}}$, then $\bm{x}(55) = \bm{x}(0)$ exactly and the dynamics are exactly periodic with period $55$. For any $s \in I$, starting from $\bm{x}(0) = \bm{x}_s$, the dynamics will converge to that trajectory, i.e.: $||\bm{x}_s (t) - \bm{x}_{s^{*}}(t)||_{\infty} \rightarrow 0$ as $t \rightarrow \infty$. This proves Theorem \ref{thm:main}.   

\begin{remark}\label{rem:perturb}
  That the parameter $s$ can belong to an interval of positive measure does not by itself imply there is an open set of initial configurations $\bm{x}(0) \in \mathbb{R}^{53}$ which fail to converge under the SBI dynamics, since for any $s$ the initial state will still consist of $5$ clusters. However, it is simple to see that, given any $\bm{x}(0)$ satisfying Theorem \ref{thm:main}, and for any sufficiently small perturbation of all initial opinions and all influence radii, the perturbed initial state will still fail to converge. Indeed, for a sufficiently small perturbation, the $5$ clusters will form at $t=1$ and, as long as they are at positions close enough to those in Table \ref{tab:one}, the clusters $C$ and $V$ will still converge to a periodic orbit of period $55$. In particular, this means we can construct non-convergent initial states consisting of $53$ agents, all with different opinions and different influence radii.   
\end{remark}

\vspace{3mm}

\section{{\bf Future work}}\label{sect:future}

Our examples, together with those in \cite{ET}, confirm that heterogeneous Hegselmann-Krause dynamics can exhibit non-convergent, periodic behaviour, whether the heterogeneity is in terms of how ``open-minded'' (different confidence radii) or how ``influential'' (different influence radii) the agents are. There are some immediate questions regarding how these recent discoveries may be refined. For example:
\\
\\
(1) How ``simple'' can a non-convergent configuration be? This paper suggests there may be no single natural way to define ``simplest''. We could try to minimize the number of agents in SBC, or the number of opinion clusters in SBI. In either model, we could try to minimize the number of distinct radii - in other words, the number of different ``types'' of agents.
\\
\\
(2) How common is non-convergent behaviour? The basic idea here would be to imagine that the initial opinion vector is chosen uniformly at random from some probability space. Of course, part of the issue here is to say what that probability space should be. However, Remark \ref{rem:perturb} suggests that for any reasonable choice of space, the probability of non-convergence will be strictly positive. Similar remarks apply to the SBC model and the example in \cite{ET}.
\\
\\
(3) Can a non-convergent trajectory be non-periodic? So far, all publicized non-convergent trajectories for both SBC and SBI are periodic.
\\
\\
(4) What, if any, are the essential differences between the SBC and SBI dynamics, at least when it comes to the issue of convergence? This is a vague question, but is motivated by the fact that, so far, research into both seems to have been motivated by an underlying assumption that they should behave similarly, in general terms. We note, however, that there is much less literature on SBI than on SBC, presumably because the original Hegselmann-Krause model was formulated as a bounded confidence, rather than bounded influence model.

\vspace{3mm}

\section{{\bf AI statement}}\label{sect:AI}

The examples in Section \ref{sect:example} were found by
\texttt{Claude Opus 5.5}. Many prompts were tried and failed until we tested the idea to look for examples consisting of differently weighted clusters, rather than single agents. The text of this paper is entirely the work of the authors.


\end{document}